\documentclass{article}
\usepackage{arxiv}
\usepackage[utf8]{inputenc}
\usepackage[T1]{fontenc}
\usepackage[ukrainian,russian,english]{babel}
\usepackage[unicode, pdftex]{hyperref}
\usepackage{url}            
\usepackage{booktabs}       
\usepackage{amsfonts}       
\usepackage{nicefrac}       
\usepackage{microtype}      
\usepackage{lipsum}		
\usepackage{graphicx}
\usepackage{amsmath}
\usepackage{amsfonts}
\usepackage{amssymb}  
\usepackage{natbib}
\usepackage{doi}
\usepackage{color}

\title{Algorithms and topological invariants for dynamic systems. II. Discrete Structures}

\author{ \href{https://orcid.org/0000-0002-7164-807X}{\includegraphics[scale=0.06]{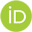}
\hspace{1mm}Alexandr O.~Prishlyak}\thanks{https://sites.google.com/view/prof-prishlyak, https://orcid.org/0000-0002-7164-807X} \\
	Department of Computer Methods of \\ 
	Mechanics and Control Processes\\
	Taras Shevchenko National University of Kyiv\\
	Kyiv, Ukraine \\
	\texttt{prishlyak@knu.ua} }

\renewcommand{\shorttitle}{Algorithms and topological invariants for dynamic systems. II. Discrete Structures}
\newtheorem{lemma}{Lemma}
\newtheorem{theorem}{Theorem}

\newtheorem{prop}{Proposition}

\newtheorem{problem}{Problem}

\hypersetup{
pdftitle={Algorithms and topological invariants for dynamic systems. II. Discrete Structures},
pdfauthor={Alexandr O. Prishlyak},
pdfkeywords={Topological classification, Morse function, Morse-Smale flow},
}

\begin{document}
\maketitle

\begin{abstract}
We construct algorithms and topological invariants that allow us to distinguish the topological type of a surface, as well as functions and vector fields for their topological equivalence.
In the first part (arXiv:2501.15657), we discused  basic concepts of diferential topology.   
In the second part we discus the main discrete topological structures used in the topological theory of dynamic systems: simplicial complexes, regular SW-complexes, Euler characteristic and homology groops, Morse-Smale complexes and handle decomposition of manifolds, Poincare rotation index of vector field, discrete Morse function and vector fields. 
\end{abstract}


\section*{Introduction}

Topological properties of manifolds, functions, and dynamical systems are often studied through the construction of topological invariants that have a discrete nature, meaning they can be described using a finite set of integers, which allows for the use of computational techniques when working with them. Morse theory enables the construction of a cell complex structure on a manifold; however, continuous mappings used to attach cells are generally difficult to encode. Therefore, triangulations are more commonly used in computer modeling. Their drawback is the large number of simplices required to construct such structures. A compromise solution between these two structures is regular cell complexes. To find all possible structures under study, efficient algorithms for their recognition are necessary. Topological invariants will be useful if there are efficient algorithms for their computation and comparison. The issues mentioned are addressed in low dimensions (up to four) in this tutorial.

The first part  (\cite{prish2025atids1}) discusses the main structures used in the topology of manifolds: vector fields, dynamical systems, Morse functions, cell decompositions, and the fundamental group.

The second part examines discrete structures for which computers can be used for their specification and manipulation.

The third part describes algorithms that allow for the recognition of manifolds and some of their properties for the structures described in the second chapter.

The fourth part is dedicated to describing possible topological structures of functions and dynamical systems on low-dimensional manifolds.


\newpage

\section{Simplicial Complexes}

\textbf{Definition.} 
A set of points $x_0, x_1,\ldots , x_n$ is called \textit{independent} if the vectors $\overrightarrow{x_0x_1}, \overrightarrow{x_0x_2},\ldots , \overrightarrow{x_0x_n}$ are linearly independent.

\textbf{Definition.}  \index{simplex}
An n-dimensional \textit{simplex} $\Delta^n$ or n-simplex is the convex hull over a system of $n+1$ independent points $(x_0, x_1,\ldots , x_n):$
$$\Delta^n=\Delta(x_0, x_1,\ldots , x_n)=\{ \sum_{i=0}^n{a_ix_i} : a_i \ge 0, \sum_{i=0}^n{a_i}=1\}.$$

The points $x_0, x_1,\ldots , x_n$ are the vertices of the simplex $\Delta(x_0, \ldots , x_n)$.

\textbf{Example.} 
0-simplices are points, 1-simplices are line segments, 2-simplices are triangles (along with their interior points), and 3-simplices are full tetrahedra.
 
\textbf{Definition.} 
Two bases of a vector space are called \textit{equivalent} if the determinant of the transition matrix from one to the other is positive. \textit{Orientation} is the class of equivalent bases.

\textbf{Definition.}  \index{orientation} The \textit{orientation} of an n-dimensional simplex is defined as the orientation of the n-dimensional plane that contains the simplex. If we fix the order of the vertices of the simplex $(x_0, x_1,\ldots , x_n)$, then the basis $\overrightarrow{x_0x_1}, \overrightarrow{x_0x_2},\ldots , \overrightarrow{x_0x_n}$ determines the orientation of the simplex. We assume that this orientation is given for $\Delta(x_0, x_1,\ldots , x_n)$. A simplex with the chosen orientation is referred to as oriented and denoted as $\delta$ or $\delta(x_0, x_1,\ldots , x_n)$. A simplex with the opposite orientation is denoted as $- \delta$.

\begin{problem}
Show that two orders of vertices define the same orientations of the simplex if they differ by an even permutation, and opposite orientations if they differ by an odd permutation.
\end{problem}

\textbf{Definition.}  \index{face} The \textit{face} of a simplex is defined as the convex hull of a subset of the vertices of the simplex. The empty set is considered a face of dimension $-1$ of any simplex. 0-dimensional faces, \index{vertex} are called \textit{vertices}, and 1-dimensional faces are referred to as \index{edge} \textit{edges}.

\textbf{Definition.}  \index{simplicial complex} A \textit{simplicial complex} $K$ (geometric simplicial complex) is defined as a set of simplices in $\mathbb{R}^m$ that satisfies the following conditions:

1) if a simplex belongs to $K$, then so does each of its faces;

2) the intersection of any two simplices in $K$ is a face of each of them;

3) every point of a simplex in $K$ has a neighborhood that intersects with a finite number of simplices.

The dimension of the complex is defined as the maximum dimension of the simplices it contains.

\textbf{Definition.}  \index{polyhedron}
The \textit{support} $|K|$ of the complex $K$ is defined as the union of all simplices contained in it. The support of the complex $K$ is also referred to as a (Euclidean) \textit{polyhedron}.

\textbf{Definition.} 
A simplicial complex is called finite if it consists of a finite set of simplices.

\begin{problem}
Prove that a simplicial complex $K$ is finite if and only if its support $|K|$ is compact.
\end{problem}

\textbf{Definition.}  \index{triangulation}
A triangulation of a topological space $X$ is a pair $(S, h)$ consisting of a simplicial complex $S$ and a homeomorphism $h: |S| \to X$. The images of simplices, vertices, and faces under the homeomorphism $h$ are also called (curvilinear) simplices, vertices, and faces. A space for which a triangulation exists is called triangulable.

\textbf{Definition.}  \index{subcomplex} \index{skeleton}
A subcomplex $L$ of a complex $K$ is a subset $L \subset K$ that is itself a complex. A subcomplex consisting of all faces of the complex $K$ whose dimension does not exceed $n$ is called the $n$-skeleton of the complex $K$ and is denoted by $K^n$.

A mapping of one simplex to another is called linear if it is a restriction of a linear mapping of the planes (spaces) containing these simplices. To define a linear mapping of a simplex, it is sufficient to specify the mapping of its vertices.

\textbf{Definition.}  \index{simplicial mapping}
Let $K$ and $L$ be simplicial complexes. A \textit{simplicial mapping} from $K$ to $L$ is a mapping $f$ of polyhedra such that the images of simplices are simplices, and the restriction of $f$ to each simplex is a linear mapping of simplices.

From the definition, it follows that in a simplicial mapping, vertices map to vertices.

\textbf{Definition.}  \index{isomorphic complexes}
Two simplicial complexes are called \textit{isomorphic} if there exists an isomorphism between their polyhedra that is a simplicial mapping, which is a homeomorphism, and the inverse mapping is also a simplicial mapping.
 
An isomorphism between simplicial complexes can be defined by a bijective mapping of vertices.

An isomorphism of a simplicial complex to itself that is not the identity is called a symmetry.

\textbf{Definition.}  \index{subdivision}
A complex $L$ is called a \textit{subdivision} of the complex $K$ if $|K|=|L|$ and each simplex of $L$ lies within some simplex of $K$.

We  divide each 1-dimensional simplex of the complex $K$ in half, each 2-simplex into 6 simplices along meridians, each 3-simplex along planes that pass through the meridian of a face and the opposite vertex, and so on. The resulting complex $K'$ is called the barycentric subdivision of the complex $K$. The second barycentric subdivision $K''$ is the barycentric subdivision for $K'$.

\begin{theorem} (about simplicial approximation). If there is a continuous mapping of polyhedra, then there exists such a subdivision $L$ of the complex $K$ and a simplicial mapping $L$ to $M$, which is homotopic to $f$.
\end{theorem}
\textbf{Definition.}  \index{algebraic simplicial complex}
An \textit{algebraic simplicial complex} on a set $V$ is defined as a collection $K \subset 2^V$ of subsets of $V$ that has the property: if $A \in K$ and $B \subset A$, then $B \in K.$
 
Each geometric simplicial complex defines an algebraic one on its vertex set. The geometric realization of an algebraic complex is an injective mapping $i: V \to \mathbb{R}^m$, which defines the geometric simplicial complex.

For finite simplicial complexes, it is convenient to take the set $V$ as the set of non-negative integers $0, 1, 2, \ldots , m$. Then one of the possible geometric realizations of the complex in $\mathbb{R}^m$ can be given by the mapping $$i(0)=(0,0,\ldots ,0), i(1)=(1,0,\ldots ,0), $$ $$i(2)=(0,1,0,\ldots ,0),\ldots ., i(m)=(0,0,\ldots , 0,1). $$

We denote by $c(K)$ the vector (sequence) consisting of the number of simplices of dimension 0, 1, 2, 3, \ldots, dim(K).

\textbf{Example.} 
Let the complex 
$$K=\{\{0\},\{1\},\{2\},\{3\},\{4\},\{5\},\{0,1\},\{0,2\},\{0, 3\},\{0,4\},\{0,5\},\{1,2\},\{2,4\},\{0,1,2\}\}.$$ Then   $$c(K)=\{6, 7,1\}.$$

If $V$ is a finite set, and $A \subset 2^V$ is some set of its subsets, then if we combine it with all the subsets of elements from $A$ (adding all the faces of the simplices), we obtain a simplicial complex, which we denote as $\left\langle A \right\rangle$ and call the complex generated by $A$. In particular, the complex from the last example can be written as:
$$
K= \left\langle \{0, 3\},\{0,4\},\{0,5\},\{2,4\},\{0,1,2\} \right\rangle.
$$
If the number of vertices does not exceed 10, we also use a shortened notation $K= \left\langle 03, 04, 05, 24, 012 \right\rangle.$

\begin{problem}
Construct a triangulation for the following surfaces:

1) torus; 2) torus with a hole; 3) torus with two holes; 4) sphere with three holes; 5) Klein bottle; 6) Möbius strip with a hole; 7) Klein bottle with two holes.
\end{problem}

\section{Regular Cell Complexes}

\textbf{Definition.}  \index{regular complex}
A cell complex $X$ is called regular if for any cell $e^k \in X$, its characteristic mapping $S^{k-1} \to X^{k-1}$ is injective, and the intersection of the closures of two cells from $X$ is either empty or the closure of another cell from $X$.

\begin{theorem}Every regular cellular complex (RCC) is triangulated.
\end{theorem}
\textbf{Proof} can be conducted by induction on the dimension of the cells. By definition, a 1-skeleton is a graph without loops and multiple edges, hence a simplicial 1-complex. Therefore, for every 2-cell, its boundary is a simplicial complex $C$. We consider this cell (along with its boundary) as a cone over $C$. Thus, we have a decomposition into simplices. Similarly, we construct decompositions into simplices for other 2-dimensional cells, then for 3-dimensional cells, and so on. As a result, we obtain a simplicial decomposition of the cellular space.

This theorem allows the RCC to be defined analogously to simplicial complexes using finite sets of elements: for each cell -- a list of cells that are part of its boundary and have dimension one less.

If the boundary of a 2-cell consists of 1-cells $\{v_1,v_2\},$ $\{v_2,v_3\},$ $ \ldots$, $\{v_k,v_1\}$, we use the abbreviation $[ v_1,v_2,v_3,\ldots ,v_k]$ for it.

As in the case of a simplicial complex, we use the abbreviated notation for the RCC, from which the full notation can be obtained by adding all cells that belong to the boundaries already listed.

\begin{figure}[ht]
\center{\includegraphics[height=4.2cm]{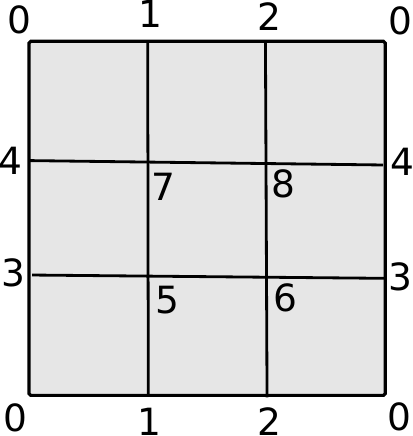}}
\caption{RCC of the 2-dimensional torus}
\label{pkkt2}
\end{figure}

\textbf{Example.} 
In figure \ref{pkkt2}, a regular cellular decomposition of the 2-dimensional torus is shown. Its abbreviated notation can be written as:
$$RCC(T^2)=\left\langle [0,1,7,4],[1,2,8,7],[2,0,4,8],[3,4,7,5], \right.$$ $$\left. [5,7,8,6], [3,6,8,4], [0,3,5,1],[1,5,6,2],[0,2,6,3] \right\rangle.$$
 
\textbf{Example.} 
For the 3-dimensional sphere, as the boundary of a 4-simplex, we obtain such an RCC $(c(RCC(S^3)=\{5, 10,10,5\})$ (the numbering of the 2-cells is shown in blue):

$$RCC(S^3)^2=\{\textcolor{blue}{1} [0,1,2],\textcolor{blue}{2}[0,1,3],\textcolor{blue}{3}[0,1,4],\textcolor{blue}{4}[0,2,3],\textcolor{blue}{5}[0,2,4], $$ $$ \textcolor{blue}{6}[0,3,4],\textcolor{blue}{7}[1,2,3], \textcolor{blue}{8} [1,2,4], \textcolor{blue}{9}[1,3,4],\textcolor{blue}{0}[2,3,4]\},$$
$$RCC(S^3)^3=\{ \{\textcolor{blue}{0,7,8,9}\},\{\textcolor{blue}{4,5,6,0}\},\{\textcolor{blue}{2,3,6,9}\},$$ $$\{\textcolor{blue}{1,3,5,8}\},\{\textcolor{blue}{1,2,4,7}\} \}.$$

Other examples of regular cell complexes are provided in the appendices.

\begin{problem}
Construct a regular cell decomposition for the following spaces:

1) immersed projective planes in three-dimensional space with one triple point and three loops, which are sets of double self-intersection points (the Boy and Gödel surfaces); 2) an immersed sphere with a set of self-intersection points -- a closed curve; 3) immersed Klein bottles with a set of self-intersection points -- a closed curve; 4) a three-dimensional torus.
\end{problem}

\section{Euler Characteristic and Homology Groups}

\textbf{Definition.}  \index{Euler characteristic}
The Euler characteristic of a compact cell complex is defined as the difference between the number of even-dimensional and odd-dimensional cells.
 
If $c(K)=\{c_0, c_1, c_2, \ldots , c_n\}$ is the characteristic vector of the complex ($c_i$ is the number of cells of dimension $i$), then the Euler characteristic is given by
$$ \chi (K) = c_0 - c_1 + c_2 - \ldots + (-1)^n c_n.$$
Since simplices can be considered as cells, for a simplicial complex, the Euler characteristic equals the number of even-dimensional simplices minus the number of odd-dimensional simplices.

\begin{theorem} \cite{hat02, prish-modtop06}
If $|K|=|L|$, then $ \chi (K) = \chi (L)$.
\end{theorem}

This theorem allows for the definition of the Euler characteristic for spaces that permit a cell decomposition (or triangulation) as the Euler characteristic of their cell decomposition.

For a simplicial complex $K$, let us denote by $C_i(K)$ the free abelian group generated by the oriented simplices of dimension $i$ from the complex $K$:
$$C_i(K) = \{ \sum_{j=1}^m a_j\Delta_j^i | a_j \in \mathbb{Z}\}.$$
We define a boundary homomorphism $\partial_i: C_i(K) \to C_{i-1}(K)$ on the generators by the formula

$$ \partial_i (\Delta (v_0,v_1,v_2, \ldots .,v_i)= \Delta (v_1,v_2, \ldots .,v_i)-\Delta (v_0,v_2,v_3, \ldots .,v_i)+$$ $$\ \ \ \ \ \ \ \ \ \ \ \ \  +\Delta (v_0,v_1,v_3, \ldots .,v_i)-\ldots  +(-1)^i  \Delta (v_0,v_1,v_2, \ldots .,v_{i-1}).$$

\begin{theorem} (Poincaré on Boundary Homomorphism) $$\partial_{i+1} \circ \partial_{i}=0.$$
\end{theorem}

\textbf{Definition.}  \index{homology group}
From this lemma, it follows that $ \text{im} \partial_{i+1} \subset \ker \partial_{i}$, and these are normal subgroups, so the quotient group is well-defined:
$$H_i(K)=\ker \partial_{i} / \text{im} \partial_{i+1}.$$
This quotient group is called the $i$-\textit{dimensional homology group of the complex} $K$.

\textbf{Definition.}  \index{chain complex}
A sequence of abelian groups $C_i$ and homomorphisms $\partial_i: C_i \to C_{i-1}, i \in \mathbb{Z}$, for which the equality $\partial_{i+1} \circ \partial_{i}=0$ holds is called a \textit{chain complex}.
 
Thus, to find the homology groups of a space, one must form a chain complex and compute its homologies. Depending on how the chain complex is constructed, different homology groups can be obtained. For instance, one can use such structures as simplicial complexes (as before), cellular spaces, Morse–Smale complexes, or discrete gradient fields. For nice spaces, such as compact simplicial complexes, all of these structures generate isomorphic homology groups in each dimension.

\begin{problem}
Find the Euler characteristic and the homology groups of the following spaces:

1) torus; 2) torus with a hole; 3) torus with two holes; 4) sphere with three holes; 5) Klein bottle; 6) Möbius strip with a hole; 7) Klein bottle with two holes.
\end{problem}

\section{Morse–Smale Complex and Handle Decomposition}
Let $M^n$ be a smooth closed manifold, and $f:M^n \to \mathbb{R}$ be a Morse function. By applying a partition of unity, one can construct a Riemannian metric that, in the neighborhood of each critical point, is induced by the Euclidean metric for the coordinate system from Morse's lemma. In this coordinate system, the gradient field  take the form
$$ \nabla f(x_1,\ldots ,x_n)=\{-2x_1,-2x_2,\ldots ,-2x_{\lambda}, 2x_{\lambda+1},\ldots ,2x_n\}.$$
Its trajectories are curves given by
$$ x_1=C_1e^{-2t}, x_2=C_2e^{-2t},\ldots ,x_{\lambda}=C_{\lambda}e^{-2t},$$ $$x_{\lambda+1}=C_{\lambda+1}e^{2t},\ldots ,x_n=C_ne^{2t}. $$
\textbf{Definition.}  \index{stable manifold}
Recall that the \textit{stable manifold} $S(p)$ (\textit{unstable manifold} $U(p)$) of a critical point $p$ is defined as the set of points that approach $p$ along their trajectories as $t\to \infty$ ($t\to -\infty$).
 
From the explicit form of the trajectories, it follows that the stable manifold $S(p)$ is determined by the condition $C_{\lambda+1}=\ldots = C_n=0$ or, equivalently, $x_{\lambda+1}= \ldots = x_{n}=0$. Thus, the index of the critical point $\lambda = \dim S(p)$. Similarly, $\dim U(p)= n-\lambda$. In particular, the stable manifold is homeomorphic to the interior of a disk of dimension $\lambda$.

\begin{figure}[ht]
\center{\includegraphics[height=4.2cm]{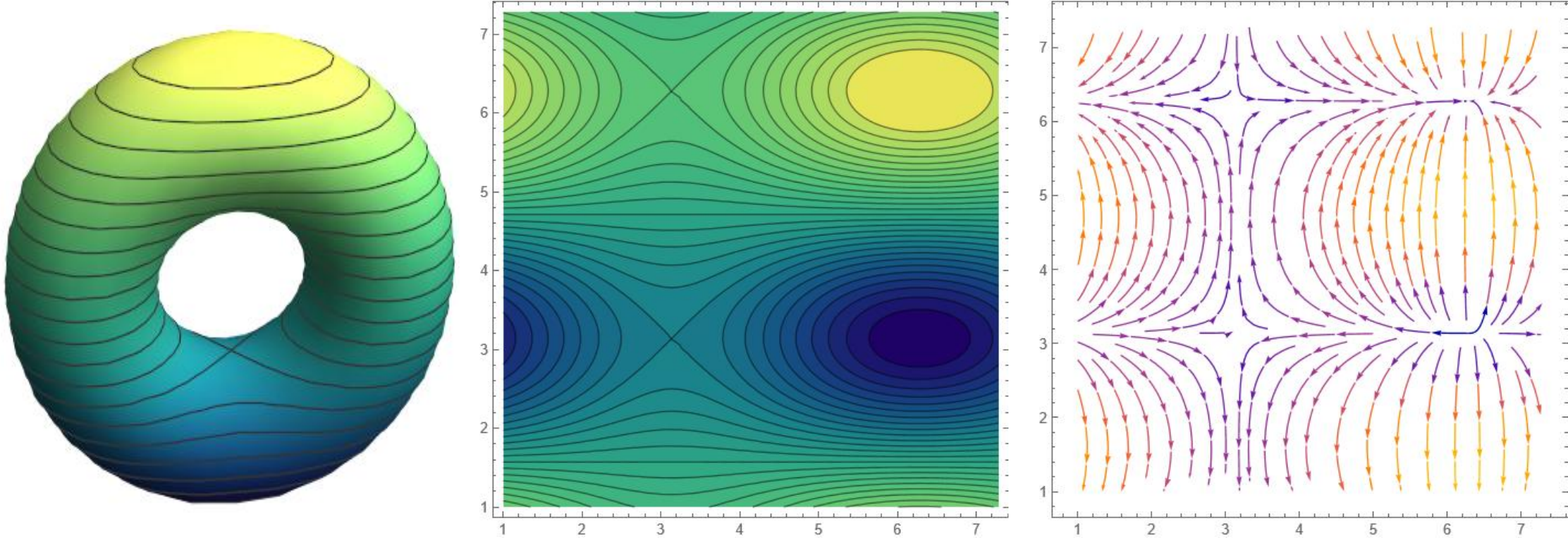}}
\caption{Gradient flow with saddle connection on the torus}
\label{h-fun-t2}
\end{figure}

\begin{figure}[ht]
\center{\includegraphics[height=4.2cm]{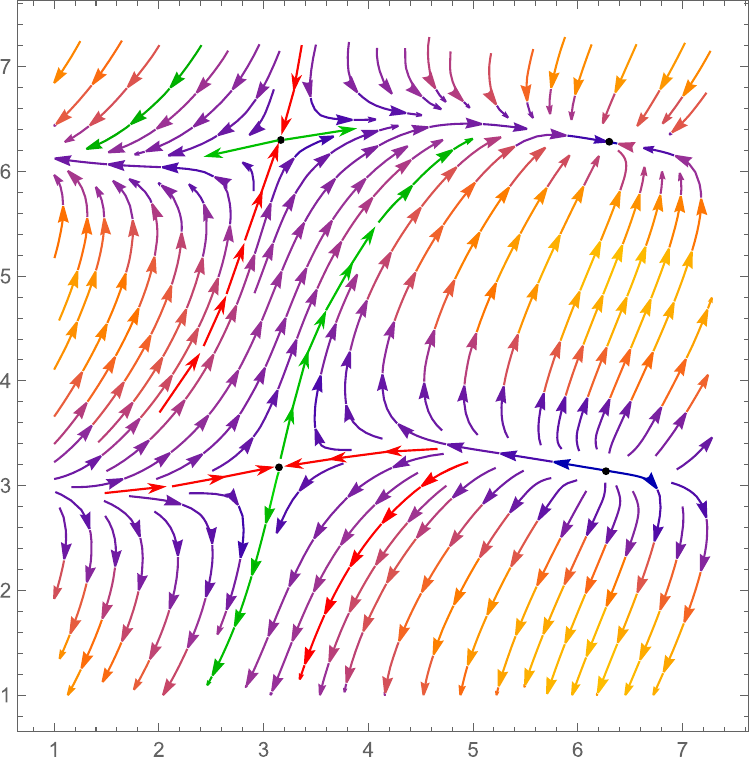}}
\caption{Bringing stable and unstable manifolds into general position}
\label{sc-bif}
\end{figure}

If necessary, by changing the Riemannian metric, one can bring stable and unstable manifolds into general position.
For example, the gradient field of the function $h(u,v)=(2+\cos u) \cos v$ on the torus in the Riemannian metric $ds^2=du^2+dv^2$ (see Fig. \ref{h-fun-t2}) has coordinates
$\{ -\sin u \cos v, -\sin v (2+\cos u)\}.$ If we slightly modify the metric, for instance to the metric $ds^2=1.25 du^2+du dv+dv^2$, we obtain a gradient field with coordinates
$$\{-\sin u \cos v - 0.5 \sin v (2+\cos u),-0.5 \sin u \cos v-1.25\sin v (2+\cos u)\}$$
and without saddle connections (see Fig. \ref{sc-bif}).

If $\text{ind} (p) =\lambda, \text{ind} (q) =\mu, \lambda \le \mu$, then in general position $S(p) \cap U(q) = \emptyset$. This means that the boundary $\partial S(p) \subset \cup_{x, \text{ind} (x)<\lambda} S(x)$.

Thus, the set of stable manifolds of critical points forms a cell decomposition of the closed manifold $M^n$. This cell complex is called the Morse–Smale complex.

Given a cell decomposition of a smooth manifold, we thicken the 0-cells (taking their regular neighborhoods), then in addition to the thickened 0-cells, we thicken the 1-cells, and in addition to the thickened 0- and 1-cells, we thicken the 2-cells, and so on. The resulting thickened cells are homeomorphic to $n$-dimensional disks and are called handles. The attachment of cells is replaced by the attachment of handles.

Let $M$ be an $n$-dimensional manifold with boundary $\partial M$.

\textbf{Definition.}  \index{handle} If $M'=M\cup H$ and there exists a homeomorphism $\varphi : D^{\lambda}\times D^{n-\lambda} \to H$ such that $\varphi(\partial D^{\lambda}\times D^{n-\lambda})=H\cap M \subset \partial M$, then $H$ is called a handle of index $\lambda$ (or a $\lambda$-handle). The mapping $\varphi$ is called the characteristic map, and the restriction $\psi =\varphi |_{\partial D^{\lambda}\times D^{n-\lambda}}$ is called the attaching map. $M'=M \cup H= M\cup _{\psi} D^{\lambda}\times D^{n-\lambda}$ is obtained from $M$ by attaching a handle of index $\lambda$. In particular, $\varphi ( D^{\lambda}\times \{0\} )$ is called the core disk or axis, $\varphi (\{0\} \times D^{n-\lambda} )$ is called the co-core disk or co-axis, $\varphi ( \partial D^{\lambda}\times \{0\} )$ is called the core sphere or a-sphere, $\varphi (\{0\} \times \partial D^{n-\lambda} )$ is called the co-core sphere or b-sphere, and $\varphi (\{0\} \times \{ 0 \} )$ is called the center of the handle.


\textbf{Definition.}  \index{handle decomposition} \textit{Handle decomposition} of a closed manifold $M$ is a sequence of inclusions $$M_0 \subset M_1 \subset M_2 \subset \ldots \subset M_m = M$$ such that $M_0$ is an $n$-dimensional disk, and $M_{i+1}$ is obtained from $M_i$ by attaching a handle. We also use the notation for handle decomposition:
$$M = H_0 \cup H_1 \cup H_2 \cup \ldots \cup H_m,$$
where $M_k = \sum_{i=0}^k H_i.$
 

For a Morse function $f$ on a smooth manifold $M$, we  construct a handle decomposition of $M$ in which each handle contains one critical point of the function $f$, and the indices of the corresponding handles and critical points are the same. Let $X_f$ be a gradient-like vector field of the function $f$ (a gradient field in some Riemannian metric on $M$), and let $p_0, p_1, \ldots, p_m$ be the critical (singular) points of the function $f$ (the field $X_f$), with $f(p_0) \le f(p_1) \le \ldots \le f(p_m)$. We set $H_0$ to be a regular neighborhood of $p_0$ such that $p_i \notin H_0$ for $i > 0$, and $\partial H_0$ is transverse to the field $X_f$; $H_i$ is a regular neighborhood ($\epsilon$-neighborhood in the corresponding metric) of the stable manifold $S(p_i)$ of the singular point $p_i$ in $M \setminus \cup_{j=0}^{i-1} H_j$, which does not contain other critical points and is transverse in $\partial H_i \setminus \cup_{j=0}^{i-1} H_j$ to the field $X_f$, for $0 < i < m$. In this constructed handle decomposition, the stable disks lie on stable manifolds, while the unstable disks lie on unstable manifolds corresponding to the singular points.

\textbf{Definition.}  Two handle decompositions of manifolds $M$ and $M'$ are called \textit{isomorphic} if there exists a homeomorphism $g: M \to M'$ that maps handles to handles, stable disks to stable disks, and unstable disks to unstable disks correspondingly.  

The handle decomposition constructed above, up to isomorphism, depends on the choice of Riemannian metric on $M$ and does not depend on the choice of neighborhoods.

\textbf{Definition.}  \index{isotopic embeddings} Two embeddings $f_0, f_1: N\to M$ are called \textit{isotopic} if there exists a homotopy $f_t$ between them such that for each $t\in [0,1]$, $f_t$ is an embedding. Two embeddings $f_0, f_1: N\to M$ are called volume isotopic if there exists a set of diffeomorphisms $g_t: M\to M, t\in [0,1]$, which depend continuously on $t$ and satisfy $g_0=\text{id}$, $f_0=g_1(f_1)$. Every volume isotopy induces a simple isotopy as a restriction to $f_0(N)$.

If handles $H_1, H_2$ are attached to $M$ via isotopic attaching maps, then $M\cup H_1$ and $M\cup H_2$ are diffeomorphic.

Using isotopies of attaching maps, one can arrange the middle and co-middle spheres of the handles into a standard position in the corresponding $\partial M_i$. In future discussions of handle decompositions, we assume that the middle and co-middle spheres intersect transversally. If there is one point of transversal intersection between these spheres, the handles are called complementary. If there are no intersections, then through isotopy, the handles can be arranged so that they do not intersect.

\textbf{Definition.}  \index{incidence index of handles} The \textit{incidence index} $\delta(H_i^{\lambda}, H_j^{\lambda -1})$ of handles $H_i^{\lambda}, H_j^{\lambda -1}$ is defined as the algebraic number of intersection points between the middle sphere of handle $H_i^{\lambda}$ and the co-middle sphere of handle $H_j^{\lambda -1}$ in $\partial M_{i-1}$.

When working with handles, the following fundamental principles are used:

\textbf{Principle 1.} Rearrangement of handles. If in a handle decomposition, neighboring handles do not intersect, they can be swapped.

\textbf{Definition.}  A handle decomposition is called \textit{proper} if the handles are attached in increasing order of indices, and handles of the same index do not intersect (they are attached simultaneously).

From any handle decomposition, a correct handle decomposition can be made by rearranging the handles. This follows from the fact that if handle $H^{\lambda}$ is attached before handle $H^{\mu}$ and $\lambda \geq \mu$, then in general position, the middle sphere $H^{\mu}$ does not intersect the core sphere $H^{\lambda}$ at the boundary of the manifold obtained after attaching $H^{\lambda}$. Therefore, by isotoping the attaching map of $H^{\mu}$, one can ensure that the handles do not intersect and change the order of attachment. Henceforth, a handle decomposition is considered correct.

\textbf{Principle 2.} Adding handles.
This operation involves sliding (isotoping) one handle along another handle of the same index. At the moment of addition, the middle and core spheres of the corresponding handles intersect at a single point.

\textbf{Principle 3.} Introducing and cancelling handles.
If the middle sphere $H^{\lambda +1}$ intersects transversally the core sphere $H^{\lambda}$ at a single point, then such handles are called complementary. Attaching a pair of complementary handles does not change the manifold (up to diffeomorphism). Therefore, if there are complementary handles in the decomposition, they can be cancelled. The reverse procedure involves introducing a pair of complementary handles.

\textbf{Construction.} (Cancelling handles of index 0). Let a connected manifold have a correct handle decomposition with more than one 0-handle. We  show that for each 0-handle, there exists a complementary 1-handle. Indeed, if for a 0-handle there were no complementary 1-handle, then this handle would be disconnected from the other 0-handles, because non-complementary 1-handles and handles of greater indices cannot connect it to the other 0-handles. Thus, all 0-handles, except for one—the last—can be cancelled.

A handle decomposition in which handles are attached in reverse order is called dual. By using the dual handle decomposition, one can cancel all but one n-handle. Therefore, on every connected closed manifold, there exists a correct handle decomposition with one 0-handle and one n-handle

\begin{prop} \textbf{(S. Smale).} One can transition from one decomposition of the manifold $M$ into handles to another using operations 1)-3) and isotopies of the middle spheres of the handles within the union of handles with lower indices.
\end{prop}

Each decomposition into handles defines a cellular decomposition of the manifold. To achieve this, we contract each handle to a middle disk. More precisely, each handle is a product of a middle and a co-middle disk, and the contraction is defined by the projection onto the first factor. These contractions provide a homotopy equivalence between the initial and final spaces. In fact, it is possible to construct a homeomorphism between them (it is sufficient to show that the regular neighborhoods of the cells, after subtracting the neighborhoods of cells of lower dimensions, are homeomorphic to the initial handles).

The constructed cellular decomposition of the manifold allows for the computation of homology groups and the fundamental group. However, this can also be done directly from the handle decomposition.

The fundamental group is computed similarly to the cellular decomposition: the generators (letters) correspond to the 1-handles, while the relations are provided by the 2-handles: each relation consists of letters corresponding to the intersection points of the middle sphere of the 2-handle with the co-middle spheres of the 1-handles, in the order dictated by the movement along the middle sphere. The degree of each letter equals the index of the corresponding intersection point.

\begin{theorem} \textbf{(Morse).} The Euler characteristic of a closed manifold equals the difference between the number of points with even and odd Morse indices of the Morse function (of the Morse-Smale complex).
\end{theorem}

\begin{problem}
Construct a handle decomposition for the following spaces:

1) torus; 2) torus with a hole; 3) torus with two holes; 4) sphere with three holes; 5) Klein bottle; 6) Möbius strip with a hole; 7) Klein bottle with two holes.
\end{problem}

\section{Poincaré–Hopf Theorem on the Indices of Vector Fields}

The degree of a mapping $f:S^1\to S^1$ is defined as the integer that equals the image of the generator (the unit) in the induced mapping $f_*:\pi_1(S^1)\to \pi_1(S^1)$ of the fundamental groups:
$$\deg f=f_*(1).$$
If a point traverses the circle once in a counterclockwise direction, the degree equals the algebraic number of turns that the image of this point makes in the counterclockwise direction (turns in the clockwise direction are counted with a negative sign). For example, the degree of a constant mapping is 0, that of the identity mapping is 1, and that of axial symmetry is $-1$. From the definition, it follows that the degrees of homotopic mappings are equal.

Let $\gamma$ be a regular closed curve in the plane, $X$ be a vector field in the plane (in some chart on a 2-manifold), and $h: \gamma \to S^1$ be a homeomorphism that preserves orientation (the direction of movement along the curve counterclockwise). We define the mapping $g: S^1 \to S^1$ by the formula:
$$g (t)= \frac{X(h^{-1}(t))}{|X(h^{-1}(t))|}.$$
Then the Poincaré rotation index of the curve $\gamma$ is defined as the degree of this mapping:
$$i_{\gamma}=\deg g.$$

\begin{figure}[ht!]
\center{\includegraphics[width=0.95
\linewidth ]{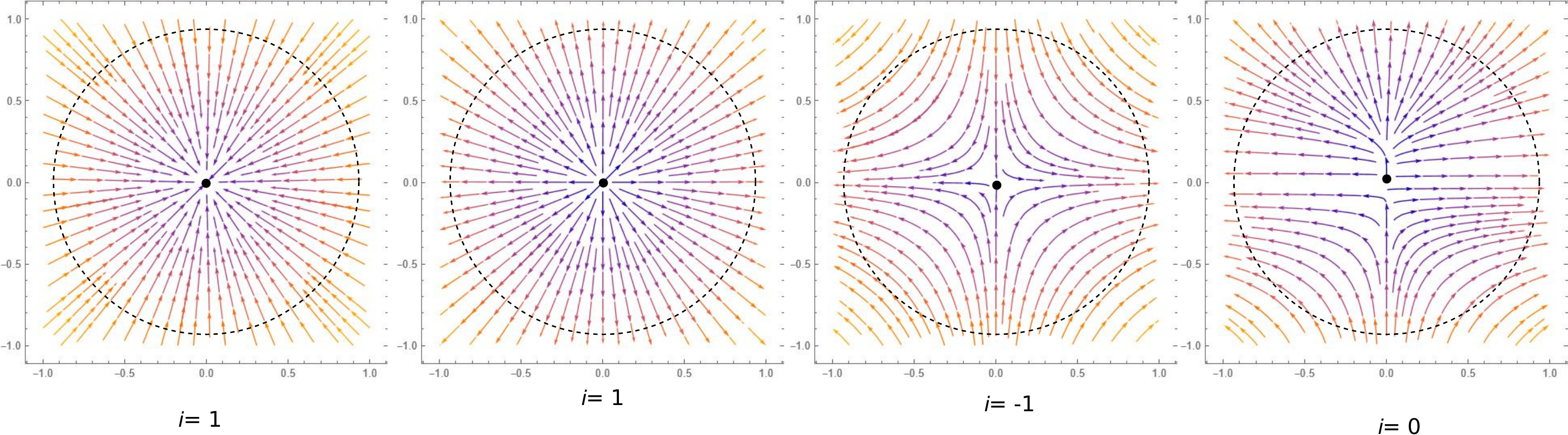}}
\caption{Examples of rotation indices}
\label{ir20}
\end{figure}

Let $p$ be an isolated singular point of the vector field $X$ in the plane (in some chart on a 2-manifold), and let $U$ be a regular neighborhood of $p$, the closure of which does not contain other singular points. The Poincaré rotation index $i_p(X)$ of the field $X$ at the point $p$ is defined as the Poincaré rotation index of the curve $\partial U$.

\textbf{Examples.} The rotation index of a saddle point is $-1$, both sources and sinks have the index of 1, and for a saddle-node, it is 0 (see Fig. \ref{ir20}). Indeed, for the source $\{x,y\}$, if a point moves counterclockwise around the unit circle, the velocity vector of the flow, which has the same coordinates as the point, also makes one complete rotation. For the sink $\{-x,-y\}$, the velocity vector of the flow is directed opposite to the radius vector of the point, and thus the corresponding symmetric point on the circle also makes one counterclockwise rotation. For the saddle $\{x,-y\}$, while the point on the unit circle traverses an angle of $\pi/2$, the velocity vector of the flow rotates by the same angle but in the opposite direction. Therefore, in total, the velocity vector makes one complete clockwise rotation. For the fourth flow $\{x,y^2\}$ (saddle-node) in the upper half-plane, the velocity vector rotates like the point and makes a rotation of $\pi$, while in the lower half-plane, it rotates in the opposite direction by the same angle. Thus, the total rotation angle equals $\pi - \pi = 0$.

\begin{theorem}
\textbf{(Poincaré–Hopf)} The sum of the Poincaré rotation indices of the  vector field on a closed manifold equals the Euler characteristic of the manifold.
\end{theorem}

\textbf{Remark}. The Poincaré-Hopf theorem is also referred to as the index theorem and the Poincaré theorem on vector fields.

To prove the theorem, we first establish several auxiliary statements.

\begin{lemma}
If a smooth vector field without singular points is defined on $S^1 = \partial D^2 \subset \mathbb{R}^2$, which is in general position with the tangent vector field and has an rotation index  0, then it can be extended to a vector field without singular points on $D^2$.
\end{lemma}

\textbf{Proof.} Let us fix one of the points of tangency of the vector field with the circle (if there were none, the winding number would be equal to 1, which contradicts the conditions of the lemma). Without loss of generality, we can assume that the tangent vector defines motion counterclockwise. We  move along the circle in this direction starting from the point of tangency and calculate the angle $\alpha$ between the vector field and the tangent vector as a continuous function. For time, we choose the polar angle from the direction to the initial point. Then we obtain the function $\alpha : [0, 2\pi] \to \mathbb{R}$, $\alpha (0) = 0$, $\alpha(2\pi)=-2\pi$ (since the winding number is equal to 0). If at 0 the function $\alpha$ is increasing, we  choose a different initial point so that it decreases (such a point exists by the intermediate value theorem). Assume that $\max \alpha = y_0>0, \alpha (t_0)=y_0$, $t_1<t_0$ is the nearest point of tangency to the left of $t_0$, and $t_2>t_0$ is the nearest point of tangency to the right of $t_0$. Then $\alpha (t_1)= \alpha (t_2)$. This means that the tangent vectors are directed in the same direction (both clockwise or both counterclockwise).
\begin{figure}[ht!]
\center{\includegraphics[width=0.6\linewidth ]{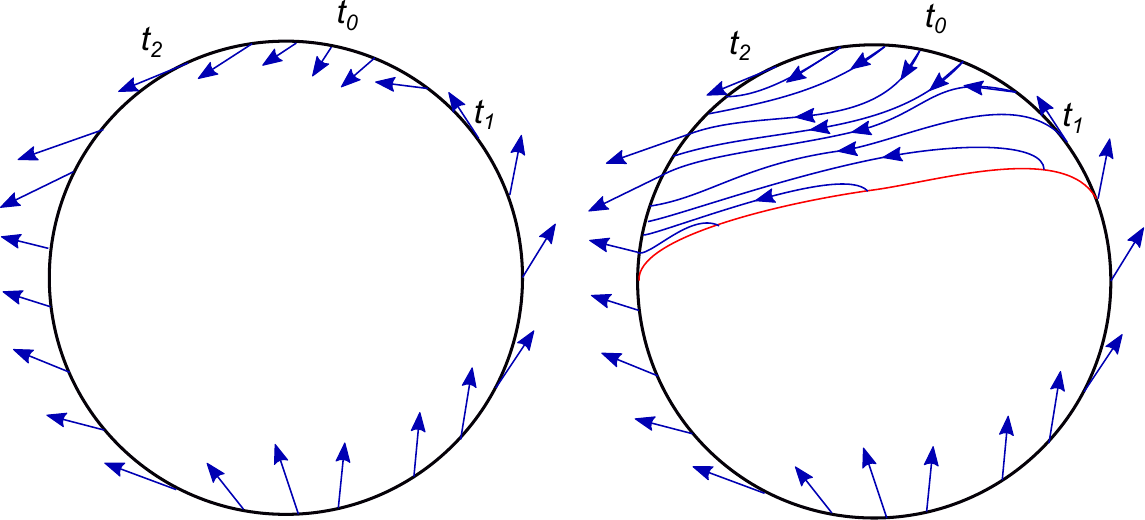}}
\caption{Reduction of two points of tangency}
\label{ph1}
\end{figure}

In Fig. \ref{ph1}, it is shown how one can extend the vector field in some neighborhood of the arc from $t_1$ to $t_2$. First, we construct a trajectory $\gamma$ that goes from $t_1$ to $t_3 > t_2$. Next, we draw trajectories from points between $t_1$ and $t_2$ so that they end between $t_2$ and $t_3$. Finally, we extend the trajectories into the neighborhood of $\gamma$ such that all of them exit at points $t > t_3$. At this point, the region left for extending the vector field is diffeomorphic to $D^2$, and the vector field on its boundary  has two fewer points of intersection. We continue this process (if necessary, also including minimum points) until only two points of intersection remain. These points split the boundary of the region into two parts—In one, the field is directed into the region, while in the other, it is directed outwards. We  fix a diffeomorphism between these arcs and connect the corresponding points with segments. The direction of movement along the segments is determined by the directions at the endpoints. After smoothing the vector field on the boundary (using a smooth homotopy), we obtain the desired vector field.

\begin{lemma}
Poincaré–Hopf theorem holds if the vector field is a Morse field.
\end{lemma}
\textbf{Proof.} For a Morse field $X$, the rotation index $i_p(X)$ is related to the Morse index $\text{ind}_p$ by the formula $$i_p(X)=(-1)^{\text{ind}_p}.$$
Then, from the formula for the Euler characteristic of the Morse–Smale cell complex, the statement of the lemma follows.
\begin{figure}[ht!]
\center{\includegraphics[width=0.4\linewidth]{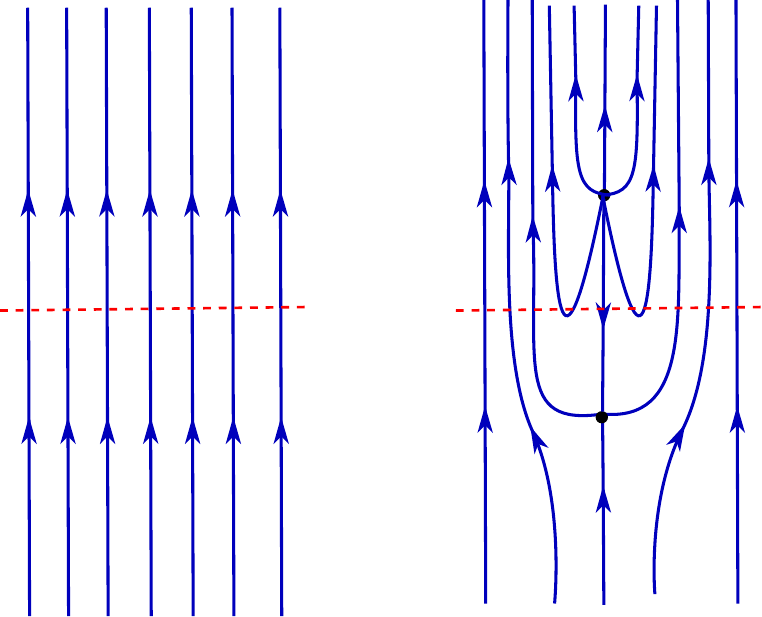}}
\caption{adding a pair of singular points: a saddle and a node}
\label{sn}
\end{figure}

Consider the operation of adding a pair of saddle and node points in a sufficiently small neighborhood of a regular point, as depicted in Fig. \ref{sn}.

\begin{lemma}
Adding a pair of saddle and node points does not change the total sum of the rotation indices.
\end{lemma}
\textbf{Proof.} The rotation index of a saddle is $-1$, while for a node, it is $1$. Therefore, their sum equals 0.

\textit{Proof of the Poincaré-Hopf Theorem.} To prove the theorem, it is sufficient to show that for any field with isolated singular points, there exists a Morse field with the same sum of the indices of rotation. For this, consider small neighborhoods of the singular points. At certain points on the boundary of the neighborhoods, we  introduce pairs of saddle points and nodes, such that one is inside the neighborhood and the other is outside, ensuring that the sum of the indices inside the neighborhood equals 0. We  replace the field inside the neighborhoods with a field without singular points as in Lemma 1. The resulting field   has only non-degenerate critical points, similar to a Morse field. We replace the obtained field with a Morse-Smale field that is close to it. Since the singular points are non-degenerate, they  remain the same in the new field. If this field has closed trajectories, we introduce a pair of saddle points and nodes on them as in Lemma 3. The resulting separatrix connections are put into general position, and we  obtain a Morse field with the same sum of the indices of rotation as the original one. The application of Lemma 2 completes the proof of the theorem.

\begin{figure}[ht!]
\center{\includegraphics[width=0.9\linewidth]{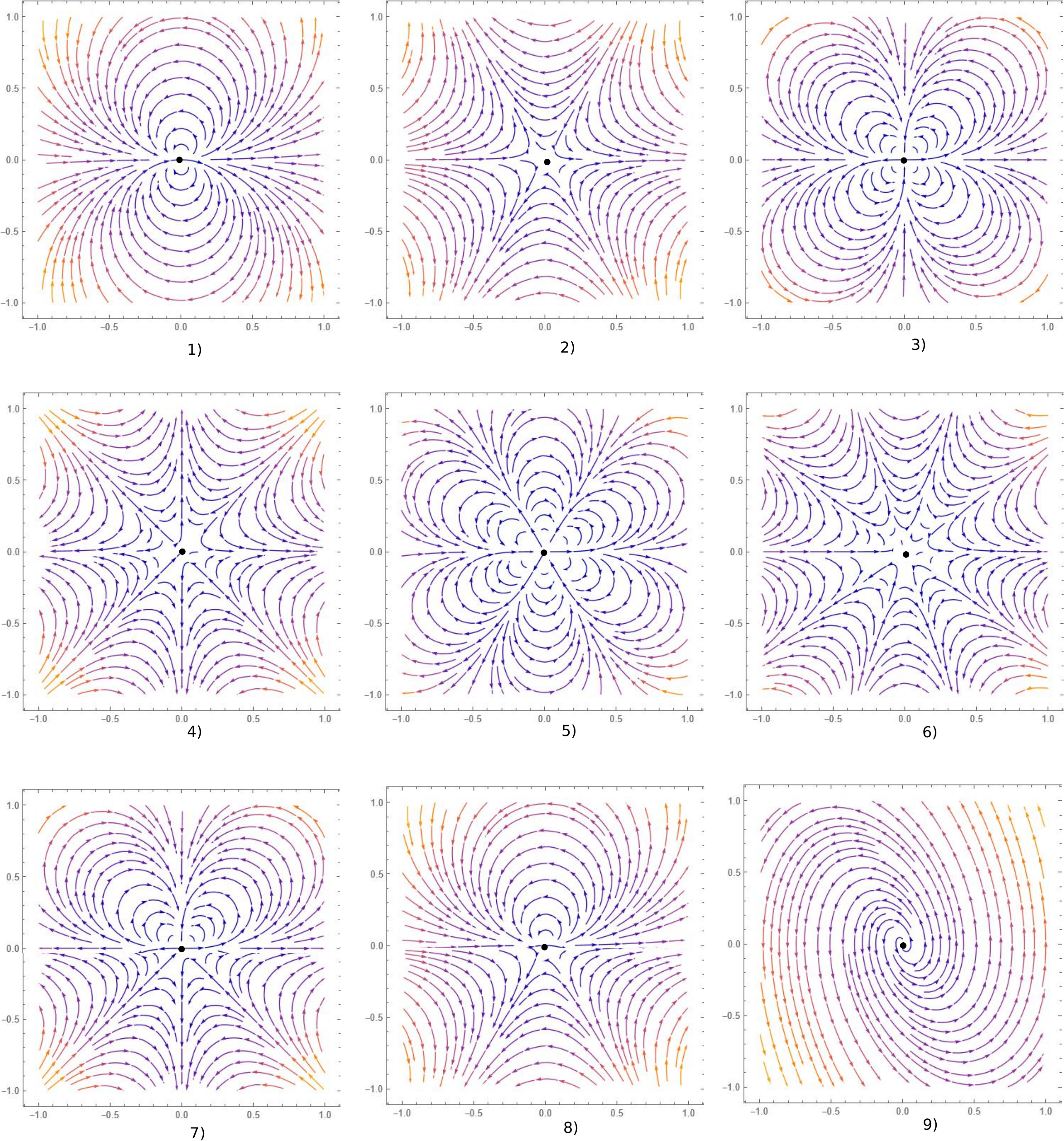}}
\caption{Singular points of vector fields in the plane}
\label{ir2}
\end{figure}

The Poincaré-Hopf theorem also holds for fields on closed n-dimensional manifolds. In this case, the index of the singular point is defined as the degree of the Gaussian map $g: S^{n-1} \to S^{n-1}$, and the degree of the mapping is defined as the algebraic number of pre-images of a regular value, that is, the difference between the number of pre-images where the mapping preserves orientation and the number of pre-images where the orientation reverses.

\begin{problem}
Find the indices of rotation for the points depicted in Figure \ref{ir2}.
\end{problem}

\section{Discrete Morse Functions and Fields}

If a cell (or simplicial) complex $K$ is given, then a discrete function is a mapping that assigns a real number to each simplex. If a manifold $M$ has a cell complex structure (triangulation), we take the value at the center of each cell (simplex) for the function $f: M \to \mathbb{R}$. In this way, we construct a discrete Morse function. Conversely, if a Morse function is given on a triangulated manifold, it determines the value of the function at the center of each simplex. We consider the first barycentric subdivision. For it, we have the function values at each vertex. By linearity, we extend the function to each simplex, thus obtaining a piecewise-linear function on the manifold. By smoothing, we can derive a smooth function from it. Similar constructions can be made for regular cell complexes.

Since smooth mappings can be approximated by piecewise-linear ones, the constructions for decomposing into handles and operations with them can be realized in the piecewise-linear category as done in \cite{rs74}.

Another approach, based on the use of discrete gradient flows, was proposed in Forman's works. This is the approach we  consider in this section, as these constructions are simpler for building algorithms and implementing them on computers.

On the set of simplices (cells) of the complex $K$, we introduce an order relation:
$$\tau < \sigma \Leftrightarrow \tau \subset \sigma, \tau \ne \sigma.$$

We denote the number of elements in the set $A$ by $|A|$.

\textbf{Definition.}  \index{discrete Morse function} A function $f: K \to \mathbb{R}$ is called a \textit{discrete Morse function} if for any simplex (cell) $\sigma$ of dimension $k$:

$$ |\tau^{k-1}: \tau^{k-1}<\sigma, f(\tau^{k-1}) \ge f(\sigma)| +
|\tau^{k+1}: \tau^{k+1}>\sigma,$$ $$f(\tau^{k+1}) \le f(\sigma)| \le 1.$$
 
\textbf{Definition.}  \index{critical simplex} \index{critical cell} A simplex (cell) $\sigma$ of dimension $k$ is called critical for the function $f: K \to \mathbb{R}$ if:

$$ |\tau^{k-1}: \tau^{k-1}<\sigma, f(\tau^{k-1}) \ge f(\sigma)| +
|\tau^{k+1}: \tau^{k+1}>\sigma,$$ $$ f(\tau^{k+1}) \le f(\sigma)| = 0.$$

Simplices that are not critical are called regular.

A discrete Morse function is called simple if it takes different values on different critical simplices (cells).

\textbf{Definition.}  \index{discrete vector field}
A discrete vector field on the complex $K$ is defined as a set of pairs $V=\{ (\sigma, \tau) \}$, where $\sigma, \tau \in K$, $\sigma <\tau$, $\dim \sigma = \dim \tau +1$, and each simplex (cell) of the complex can belong to no more than one pair. Each such pair is called a vector and is represented by an arrow directed from the center of the first simplex (cell) to the center of the second.

\textbf{Definition.}  \index{gradient field}
The gradient vector field $V_f$ of a discrete Morse function is defined as the field $$V_f=\{ \sigma^{(k)}, \tau^{(k+1)}: \sigma^{(k)}<\tau^{(k+1)}, f(\sigma^{(k)}) \ge f( \tau^{(k+1)}).$$

\begin{figure}[ht!]
\begin{center}
\begin{minipage}[h]{0.28\linewidth}
\center{\includegraphics[width=1\linewidth ]{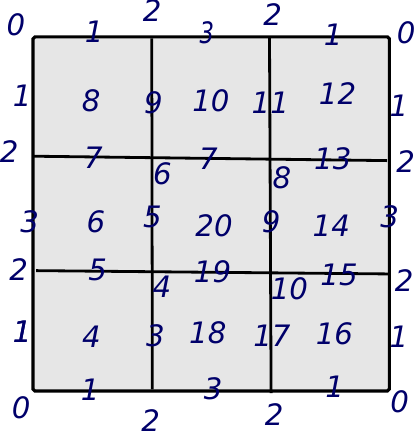}} a) \\
\end{minipage}
\  \ \ \ \ \ \ \ \ \ \ \  \ \ \ 
\begin{minipage}[h]{0.26\linewidth}
\center{\includegraphics[width=1\linewidth]{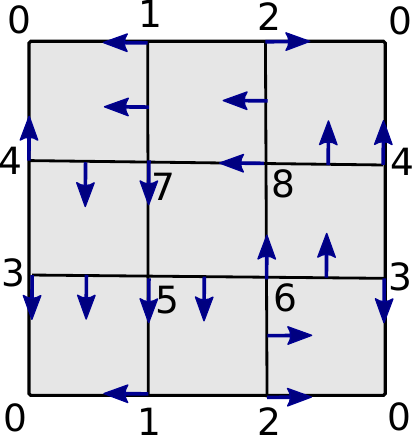}} \\b)
\end{minipage}
\end{center}
\caption{The discrete Morse function a) and its gradient field b) on the torus}
\label{dmft2}
\end{figure}

Figure \ref{dmft2} shows the discrete Morse function on the torus and its gradient field. The critical cells are: 1) vertex 0, 2) edges $\{1,2\}$ and $\{3,4\}$, 3) 2-cell $[5,6,8,7]$.

The index of a critical cell (simplex) is its dimension.

\textbf{Definition.}  \index{gradient path}
The \textit{gradient path} of a discrete vector field $V$ is defined as a sequence of cells (simplexes)
$$ \tau_0, \sigma_1, \tau_1, \ldots, \sigma_{m-1}, \tau_{m-1},\sigma_m$$
such that $$(\sigma_i, \tau_i) \in V, \sigma_{i} < \tau_{i-1}, 1 \le i \le m, \dim \sigma_0 = \dim \sigma_m.$$

Examples of gradient paths in \ref{dmft2} b): 1) two paths from $\{1,2\}$ to $\{0\}$ -- $\{1,2\},\{1\},\{0,1\},\{0\}$ and $\{1,2\},\{2\},\{0,2\},\{0\}$; 2) two paths from $[5,6,8,7]$ to $\{1,2\}$ -- $[5,6,8,7],\{5,6\},[1,2,6,5],\{1,2\}$ and $[5,6,8,7],\{5,6\},[1,2,6,5],\{2,6\}$, \ $[2,0,3,6],\{3,6\},[6,3,4,8],$ $\{8,4\},[8,4,0,2],\{2,8\},[8,2,1,7],\{2,1\}$.

By consistently aligning orientations on adjacent cells, we find that the orientation of the initial cell induces the orientation of the last cell in the path.

To compute the homology groups, we form a chain complex where the group of $k$-dimensional chains is generated by critical cells (simplexes) of dimension $k$. The boundary homeomorphisms are defined by the incidence indices between two cells of adjacent dimensions, which equals the algebraic number of gradient paths from one critical cell to another. In the previous examples, the two paths from $\{1,2\}$ to $\{0\}$ induced different orientations; therefore, the incidence index between $\{1,2\}$ and $\{0\}$ is zero. Similarly, it is zero for the pair of cells $[5,6,8,7]$ and $\{1,2\}$. The incidence indices  also are zero between $\{3,4\}$ and $\{0\}$, as well as between $[5,6,8,7]$ and $\{3,4\}$.
Thus, we have the following chain groups:
$$C_0=\mathbb{Z}, \ \ C_1=\mathbb{Z}\oplus\mathbb{Z},\ \ C_2=\mathbb{Z}.$$

All boundary homomorphisms are zero, thus $\ker \partial_i = C_i$, $\partial_i(C_i) = 0$, $H_i = C_i$. Therefore, for the torus

$$H_0(T^2) = \mathbb{Z}, \ \ H_1(T^2) = \mathbb{Z} \oplus \mathbb{Z}, \ \ H_2(T^2) = \mathbb{Z}.$$

\begin{problem}
Construct a discrete Morse function and a discrete gradient field for regular cell decompositions of such spaces with one critical 0-cell and at most one critical 2-cell:

1) torus; 2) torus with a hole; 3) torus with two holes; 4) sphere with three holes; 5) Klein bottle; 6) Möbius strip with a hole; 7) Klein bottle with two holes.
\end{problem}

\begin{problem}
For the gradient fields constructed in the previous section, find the gradient paths that start or end at critical 1-cells and use them to find the homology groups of the corresponding spaces.
\end{problem}

\section{Literature Review}







Simplicial complexes, Euler characteristics, and homology groups are fundamental concepts in homology theory courses \cite{hat02, prish-modtop06}. Regular cell complexes and discrete Morse functions are discussed in discrete Morse theory \cite{for98}. 








Topological properties of functions and vector fields have been studied in many works. Among them, we would like to highlight the works of the Kyiv topological school: 
textbooks \cite{prish2002Morse, prish-modtop06, prish2015top, prish23algkom}, 
scientific articles by the author on function topology \cite{Prishlyak1993, prishlyak1998, prishlyak1999equivalence, prishlyak2000conjugacy, prishlyak2001conjugacy, prish2001top, prishlyak2002topological1, prishlyak2002morse, prishlyak2003regular, prishter2024}, vector fields \cite{Prish1997vec, Prishlyak2001, prishlyak2002morse1, Prishlyak2002, prishlyak2003topological, prishlyak2003sum, prishlyak2005complete, Prishlyak2007}, and other geometric objects \cite{Prishlyak1994, prishlyak1997graphs, Prishlyak1999}, as well as the works of his students: K. Myshchenko \cite{prishlyak2007classification}; N. Budnytska \cite{Bud2008knu}; D. Lychak \cite{lychak2009morse}; A. Bondarenko \cite{Bond2012mfat}; O. Vyatychaninova \cite{VyatP2013Mol}; Bohdana Hladysh \cite{Hladysh2016, hladysh2017topology,
Hladysh2019, PLH2023}; A. Prus \cite{Prishlyak2017, prishlyak2020three, Prishlyak2021}; V. Kiosak \cite{KPL2022}; S. Bilun \cite{bilun2002closed, PBP23}
V. Lisikevich \cite{LisP2013KNU}, I. Ivaniuk \cite{IvanPrish2014-func-deform3, IP2015knu}, D. Skotchko \cite{
PS2016-F-atoms}, M. Loseva \cite{Losieva2017, Prishlyak2019, Prishlyak2020}, Z. Kybalko \cite{Kybalko2018}, K. Khatamian \cite{Hatamian2020}.

\bibliographystyle{unsrtnat}
\bibliography{prish}
\end{document}